\documentclass[12pt]{amsart}
\usepackage{latexsym, amssymb, amsmath, amscd}
 \usepackage{eucal}

\theoremstyle{plain}
    \newtheorem{rema}{Remark}[section]
    \newtheorem{propo}[rema]{Proposition}

   \newtheorem{theo}[rema]{Theorem}
 
   \newtheorem{defi}[rema]{Definition}
    \newtheorem{lemma}[rema]{Lemma}
    
     \newtheorem{exam}[rema]{Example}
  
\newtheorem{quest}[rema]{Question}

	\newcommand{\p}{\partial}

 \newcommand{\pf}{{\it Proof:}\hspace{2ex}}
 
 \newcommand{\epfv}{\hspace{1em}$\Box$\vspace{1em}}

\newcommand{\bZ}{{\mathbb Z}}
\newcommand{\bR}{{\mathbb R}}
\newcommand{\bQ}{{\mathbb Q}}
\newcommand{\bN}{{\mathbb N}}

\newcommand{\poly}{{\rm{Poly\,}}}
\newcommand{\supp}{{\rm {Supp\,}}}

\newcommand{\dive}{{\rm{div}\,}}

\newcommand{\bC}{\mathbb C}

\newcommand{\im}{{\operatorname{Im}\,}}

\renewcommand{\theequation}{\thesection.\arabic{equation}}
\renewcommand{\therema}{\thesection.\arabic{rema}}
\title[Images of Locally Finite Derivations]
{Images of Locally Finite Derivations of Polynomial 
Algebras in Two Variables}

\author{Arno van den Essen, David Wright and Wenhua Zhao}      
   \date{\today}

\address{A. van den Essen, Department of Mathematics, Radboud University Nijmegen,  Postbus 9010,   
6500 GL Nijmegen, The Netherlands. 
E-mail: essen@math.ru.nl }

\address{
D. Wright, Department of Mathematics, Washington University in St. Louis, St. Louis, MO. 63130, USA. \, 
E-mail: wright@math.wustl.edu}
\address{
W. Zhao, Department of Mathematics, 
Illinois State University, Normal, IL 61790-4520, USA. \, 
E-mail: wzhao@ilstu.edu.}

\begin{document}

\begin{abstract}
In this paper we show that the image of any locally finite $k$-derivation of the polynomial algebra $k[x, y]$ in two variables over a field $k$ of characteristic zero is a Mathieu subspace. We also show that the two-dimensional Jacobian conjecture is equivalent to the statement that the image $\im D$ of every $k$-derivation $D$ of $k[x, y]$ such that $1\in \im D$ and $\dive D=0$ 
is a Mathieu subspace of $k[x, y]$.
\end{abstract}

   

\thanks{The third-named author has been partially supported 
by the NSA Grant H98230-10-1-0168}

 \bibliographystyle{alpha}
    \maketitle


\renewcommand{\theequation}{\thesection.\arabic{equation}}
\renewcommand{\therema}{\thesection.\arabic{rema}}
\setcounter{equation}{0}
\setcounter{rema}{0}
\setcounter{section}{0}

\section{\bf Introduction}\label{S1}

Kernels of derivations have been studied in many papers. On the other hand, only a few results are known concerning images of derivations. 

In this paper we consider the question if the image of a derivation of a polynomial algebra in two variables over a field $k$ is a Mathieu subspace of the polynomial algebra.

The notion of the Mathieu subspaces was introduced recently 
by the third-named author in \cite{GIC} in order to study 
the Mathieu conjecture \cite{M}, the image 
conjecture \cite{IC} and the Jacobian conjecture 
(see \cite{BCW} and \cite{E1}). 
We will recall its definition in Section 2 below.

Throughout this paper we fix the following notation:
$k$ is a field of characteristic zero and $x, y$ 
are two free commutative variables. We denote by 
$A$ the polynomial algebra $k[x,y]$
over the field $k$.

The contents of the paper are arranged as follows.

In Section \ref{S2} we recall some facts concerning Mathieu subspaces 
and show that the image of a $k$-derivation of $A$ needs not be a Mathieu subspace (see Example \ref{Exam2.4}). 

In Section \ref{S3} we prove in Theorem \ref{MainThm-1} that for 
every locally finite $k$-derivation $D$ of $A$, the image $\im D$ is a Mathieu subspace. Finally in Section \ref{S4} we show in Theorem \ref{MainThm2} that the 
two-dimensional Jacobian conjecture is equivalent to the following: {\it if $D$ is a $k$-derivation of $A$ with $\dive D=0$ such that $1\in \im D$, 
then $\im D$ is a Mathieu subspace of $A$. }

\renewcommand{\theequation}{\thesection.\arabic{equation}}
\renewcommand{\therema}{\thesection.\arabic{rema}}
\setcounter{equation}{0}
\setcounter{rema}{0}

\section{\bf Preliminaries}\label{S2}

We start with the following notion introduced 
in \cite{GIC}.

\begin{defi}\label{Def-MS}
Let $R$ be any commutative $k$-algebra and 
$M$ a $k$-subspace of $R$. Then $M$ is a Mathieu subspace of $R$ if the following condition holds: if $a\in R$ is such that $a^m\in M$ for all $m\ge 1$, then for any $b\in R$, there exists an $N\in \bN$ such that $ba^m\in M$ for all $m\ge N$.
\end{defi}

Obviously every ideal of $R$ is a Mathieu subspace of $R$. However not every Mathieu subspace of $R$ is an ideal of $R$. Before we give some examples, we first recall the following simple lemma proved in Lemma $4.5$, \cite{GIC}, 
which will be very useful for our later arguments. 
For the sake of completeness, we here also include a proof.

\begin{lemma}\label{OneLemma}
If $M$ is a Mathieu subspace  of $R$ and $1\in M$, 
then $M=R$.
\end{lemma} 
\pf
Since $1\in M$, it follows that $1^m=1\in M$ for all $m\ge 1$.  Then for every $a\in R$, $a=a 1^m\in M$ for all large $m$. Hence $R\subseteq M$ and $R=M$.
\epfv

\begin{exam}
Let $R\!:=k[t, t^{-1}]$ be the algebra of Laurent polynomials 
in the variable $t$. For each $c\in k$, 
let $D_c$ be the differential operator $\frac{d}{dt}+ct^{-1}$ of $R$.  
Then  $\im D_c\!:=D_c R$ is a Mathieu subspace 
of $R$ if and only if 
$c\not \in \bZ$ or $c=-1$.   
\end{exam}

Note that the conclusion above follows directly by applying Lefschetz's principle to Proposition $2.6$ \cite{GIC}. Since Proposition $2.6$ in \cite{GIC} is for multi-variable case and its proof is quite involved, we here include a self-contained proof for the one variable case.

\pf Note first that for any $m\in \bZ$, 
$D_c t^m=(m+c)t^{m-1}$. So, if 
$c \not \in \bZ$, then $\im D_c=R$. 
Hence a Mathieu subspace of $R$.

If $c \in\bZ$ but $c\ne -1$, then $D_c t=(1+c)\ne 0$. 
So $1\in \im D_c$. Since $D_c t^{-c}=(-c+c)t^{-c-1}=0$, 
it is easy to see that $t^{-c-1}\not \in \im D_c$. 
Hence $\im D_c\ne R$. Then by Lemma \ref{OneLemma}, 
$\im D_c$ is not a Mathieu subspace of $R$.

Finally, assume $c=-1$. Since 
$D_{-1}t^m=(m-1)t^{m-1}$ for all $m\in \bZ$, 
it is easy to see that $\im D_{-1}$ is 
the subspace of the Laurent polynomials in $R$ 
without constant term. Then 
by the Duistermaat-van der Kallen 
theorem \cite{DK}, 
$M$ is a Mathieu subspace of $R$.
\epfv

Note that when $c=-1$, $\im D_{-1}$ is a 
Mathieu subspace of $R$. But it clearly is not an ideal 
of $R$. For more examples of Mathieu subspaces 
which are not ideals, see Section $4$ 
in \cite{GIC}.

When $c=0$, we see that $\im d/dt$ is not a Mathieu subspace  of $R$. Now observe that $k[t, t^{-1}]\simeq k[x, y]/(xy-1)$, where $t$ corresponds to the class of $x$ and $t^{-1}$ to the class of $y$. Then the derivation $d/dt$ of $R$ can be lifted to a $k$-derivation $D$ of $k[x, y]$, which maps $x$ to 
$\frac{d}{dt} \, t=1$ and $y$ to $\frac{d}{dt}t^{-1}=-t^{-2}$, i.e., $-y^2$. This leads to the following example. 

\begin{exam}\label{Exam2.4}
Let $D=\p_x-y^2\p_y$. Then $\im D$ is not a 
Mathieu subspace  of $k[x, y]$.
\end{exam} 
\pf Note that $1=Dx\in \im D$. However $y\not \in \im D$ 
since for any $g\in k[x, y]$ the $y$-degree of $Dg$ 
can not be $1$. So by Lemma \ref{OneLemma}, 
$\im D$ is not a Mathieu subspace  of $k[x, y]$.
\epfv

The following lemma will also be needed 
in Section \ref{S3}. 

\begin{lemma}\label{Lef-Lemma}
Let $R$ be any $k$-algebra, $L$ a field extension 
of $k$ and $M$ a $k$-subspace of $R$. Assume that 
$L\otimes_k M$ is a Mathieu subspace of the 
$L$-algebra $L\otimes_k R$. Then 
$M$ is a Mathieu subspace of 
the $k$-algebra $R$.
\end{lemma}
\pf We view $L\otimes_k R$ as a $k$-algebra in the obvious way. 
Since $L\otimes_k M$ is a Mathieu subspace of the $L$-algebra 
$L\otimes_k R$, from Definition \ref{Def-MS} it is easy to see that 
$L\otimes_k M$ (as a $k$-subspace) is also a Mathieu subspace 
of the $k$-algebra $L\otimes_k R$. 

Now we identify $R$ with the $k$-subalgebra $1\otimes_k R$ 
of the $k$-algebra $L\otimes_k R$. Then from Definition 
\ref{Def-MS} again, it is easy to check that 
the intersection $(L\otimes_k M)\cap R=M$ 
is a Mathieu subspace of $R$.
\epfv

Note that by the lemma above, when we prove that a 
$k$-subspace of a polynomial algebra over $k$ is a Mathieu subspace of the polynomial algebra, 
we may freely replace $k$ by any field extension of $k$. 
For instance, we may assume that $k$ is algebraically 
closed. 

To conclude this section we recall a result from  \cite{EWZ} which will be used in Section \ref{S3} below.

Let $z=(z_1, z_2, ..., z_n)$ be $n$ commutative 
free variables and $k[z, z^{-1}]$ the algebra of Laurent 
polynomials in $z_i$ $(1\le i\le n)$. 
For any non-zero 
$f(z)=\sum_{\alpha\in \bZ^n}c_\alpha z^\alpha \in 
k[z, z^{-1}]$, we denote by $\supp (f)$ the {\it support} 
of $f(z)$, i.e., the set of all $\alpha\in \bZ^n$ such that 
$c_\alpha\ne 0$, and $\poly (f)$ the ({\it Newton}) 
{\it polytope} of $f(z)$, i.e., the convex hull of 
$\supp (f)$ in $\bR^n$.  

\begin{theo}\label{DensityThm} $($\cite{EWZ}$)$ 
Let $0\ne f\in k[z, z^{-1}]$ and $u$ any rational point, 
i.e., a point with all coordinates being rational, 
of $\poly (f)$. 
Then there exists $m\ge 1$ such that 
$(\bR_+ u) \cap \supp (f^m) \ne \emptyset$.
\end{theo}

\renewcommand{\theequation}{\thesection.\arabic{equation}}
\renewcommand{\therema}{\thesection.\arabic{rema}}
\setcounter{equation}{0}
\setcounter{rema}{0}

\section{\bf Images of Locally Finite Derivations of $k[x, y]$}\label{S3}

Let $D$ be any $k$-derivation of $A(=k[x, y])$. Then $D$ is said to be {\it locally finite} if for every $a\in A$ 
the $k$-vector space spanned by the elements 
$D^ia$ $(i\ge 1)$ is finite dimensional.  

The main result of this section is the following theorem.

\begin{theo}\label{MainThm-1}
Let $D$ be any locally finite $k$-derivation of $A$.
Then $\im D$ is a Mathieu subspace of $A$.
\end{theo}

To prove this theorem, we need the following result, 
which is Corollary 4.7 
in \cite{E2}. 

\begin{propo}\label{LFDer}
Let $D$ be any locally finite $k$-derivation of $A$. 
Then up to the conjugation by a $k$-automorphism of $A$, 
$D$ has one of the following forms:  
\begin{enumerate}
  \item[i)] $D=(ax+by)\p_x+(cx+dy)\p_y$ 
 for some $a, b, c, d\in k$;
\item[ii)] $D=\p_x+by\p_y$ for some $b\in k$;
\item[iii)] $D=ax\p_x+(x^m+amy)\p_y$ for some $a\in k$ and $m\ge 1$;
\item[iv)] $D=f(x)\p_y$ for some $f(x)\in k[x]$.
\end{enumerate}  
\end{propo}

\begin{lemma}\label{L-3.3}
With the same notations as in Proposition \ref{LFDer}, the following statements hold. 
\begin{enumerate}
  \item[(a)] If $D$ is of type $ii)$, then $D$ is surjective. 
\item[(b)] If $D$ is of type $iii)$,  then 
\begin{align}
\im D=
\begin{cases}
(x^m) &\mbox{if $a= 0$.}\\
(x, y) &\mbox{if $a\ne 0$.}  
\end{cases}
\end{align}

\item[(c)] If $D$ is of type $iv)$, then $\im D=(f(x))$.
\end{enumerate}  
\end{lemma}
\pf 
$(a)$ is well-known, see \cite{C} or \cite{F} (p.\,96).
$(c)$ is obvious, so it remains to prove $(b)$.

If $a=0$, then $D=x^m\p_y$, and hence $\im D=(x^m)$. 
So assume $a\ne 0$. Replacing $D$ by $a^{-1}D$ (without changing  
the image $\im D$), we may assume that $D=(x\p_x+my\p_y)+b x^m\p_y$ 
for some nonzero $b\in k$. Observe that for any $i, j\in \bN$, 
we have 
\begin{align}\label{L-3.3-pe1}
D(x^iy^j)=(i+mj)x^iy^j+jbx^{m+i}y^{j-1}.
\end{align}

Next we use induction on $j\ge 0$ to show that 
$x^iy^j\in \im D$ whenever $i+j>0$.

First, assume $j=0$. Then by Eq.\,(\ref{L-3.3-pe1}), 
we have $Dx^i=ix^i$, and hence 
$x^i\in \im D$ for all $i\ge 1$. 

Now assume $j\ge 1$. Since $m\ge 1$, we have $m+i\ge 1$ for all $i\ge 0$.  
Then by the induction assumption, $jbx^{m+i}y^{j-1}\in \im D$ 
for all $i\ge 0$. Combining this fact with Eq.\,(\ref{L-3.3-pe1}), 
we get $x^iy^j\in \im D$ since $i+mj\ne 0$ for all $i\ge 0$. 
Hence we have proved that $x^iy^j\in \im D$ if $i+j>0$. 
Note that $1$ does not lie in 
$\im D$ since this space is contained in the ideal generated by
$x$ and $y$. Therefore we have $\im D=(x, y)$.
\epfv
%
%
%
%

\begin{lemma}\label{L-3.4}
Let $z=(z_1, z_2, ..., z_n)$ be $n$ free commutative variables and $D\!:=\sum_{i=1}^n a_iz_i\p_{z_i}$ for some $a_i\in k$ 
$(1\le i\le n)$. Then $\im D$ is a Mathieu subspace of $k[z]$.
\end{lemma}

Note that $D$ in the lemma is a locally finite 
derivation of the polynomial algebra $k[z]$. To show the lemma, let's first recall the following well-known results.

\begin{lemma}\label{ExtraLemma} 
For any polynomials $f, g\in k[z]$ and a 
positive integer $m\ge 1$, we have 
\begin{align} 
\poly(fg)&=\poly(f)+\poly(g), 
\label{ExtraLemma-e1} \\
\poly(f^m)&=m\poly(f),  
\label{ExtraLemma-e2} 
\end{align} 
where the sum in the first equation above  denotes the Minkowski sum of polytopes.
\end{lemma}
\pf Eq.\,(\ref{ExtraLemma-e1}) is well-known, which was first proved by A. M. Ostrowski \cite{O1} in 1921 (see also  
Theorem VI, p.\,226 in \cite{O2} or Lemma $2.2$, p.\,11 in \cite{Stu}). To show Eq.\,(\ref{ExtraLemma-e2}), one can first check easily that the polytope $m\poly(f)$ and the polytope obtained by taking the Minkowski sum of $m$ copies of $\poly(f)$ actually share the same set of extremal vertices, namely, the set of the vertices $m v_i$, where 
$v_i$ runs through all extremal vertices of $\poly(f)$. Consequently, these two polytopes coincide. 
Then from this fact and 
Eq.\,(\ref{ExtraLemma-e1}), we see that 
Eq.\,(\ref{ExtraLemma-e2}) follows.     
\epfv

\underline{\it Proof of Lemma \ref{L-3.4}:} If all $a_i$'s are zero, then $D=0$ and $\im D=0$. Hence the lemma holds in this case.  So, we assume that 
not all $a_i$'s are zero. 

Let $S$ be the set of integral 
solutions $\beta\in \bZ^n$ of 
the linear equation $\sum_{i=1}^n a_i\beta_i=0$.  
Note that $S\ne \emptyset$ (since $0\in S$) and 
is a finitely generated $\bZ$-module. 
Let $V$ be the subspace of $\bR^n$ 
spanned by elements of $S$ over $\bR$. 
Then $V$ is a $\bR$-subspace 
of $\bR^n$ with $r\!:=\dim_\bR V<n$. 
Furthermore, $V$ can be described 
as the set of common solutions 
of some linear equations with 
rational coefficients, since clearly the $\bQ$-vector space
generated by the $\bZ$-generators of $S$ can.


Note also that for any $\beta=(\beta_1, \beta_2, ..., \beta_n)
\in \bN^n$, we have $D z^\beta=(\sum_{i=1}^n a_i\beta_i) z^\beta$. Hence, for any $\beta\in \bN^n$, the monomial 
$z^\beta \in \im D$ iff $\beta \not \in S$, or equivalently, 
$\beta\not \in V$. 
Consequently, for any $0\ne h(z)\in \bC[z]$, we have
\begin{align}\label{L-3.4-pe1}
h(z)\in \im D \Leftrightarrow \supp (h)\cap V=\emptyset.
\end{align}

Now, let $0\ne f(z)\in \bC[z]$ such that $f^m\in \im D$ for all 
$m\ge 1$. We claim $\poly (f) \cap V=\emptyset$.
 
Assume otherwise. Since all vertices of 
the polytope $\poly (f)$ are rational 
(actually integral), every face of $\poly (f)$ 
can be described as the set of common solutions 
of some linear equations with rational 
coefficients. Since this is also the case 
for $V$ (as pointed above) and 
$\poly (f) \cap V\ne\emptyset$ (by our assumption), 
it is easy to see that there exists at least one 
rational point $u\in \poly(f)\cap V$. 
Then by Theorem \ref{DensityThm}, 
there exists $m\ge 1$ such that 
$(\bR_+ u) \cap \supp (f^m)\ne \emptyset$, and 
by Eq.\,(\ref{L-3.4-pe1}), $f^m\not \in \im D$. 
Hence, we get a contradiction. 
Therefore, the claim holds.

Finally, we show that $\im D$ is a Mathieu subspace as follows. 

Let $f(z)$ be as above and $d$ the distance between $V$ and 
$\poly(f)$. Then by the claim above and the fact that 
$\poly(f)$ is a compact subset of $\bR^n$, we have $d>0$. Furthermore, for any $m\ge 1$, by Eq.\,(\ref{ExtraLemma-e2}) we have $\poly(f^m)=m\poly(f)$. Hence, the distance between $V$ and $\poly(f^m)$ is given by $dm$. 

Now let $h(z)$ be an arbitrary 
element of $k[z]$. Note that 
by Eqs.\,(\ref{ExtraLemma-e1})  
and (\ref{ExtraLemma-e2}) we have 
$\poly(f^m h)=m\poly(f)+\poly(h)$ 
for all $m\ge 1$. Hence, for 
large enough $m$, the distance between $V$ and 
$\poly(f^m h)$ is positive, whence  
$\poly (f^mh)\cap V=\emptyset$. 
In particular, $\supp (f^mh)\cap V=\emptyset$, 
and by Eq.\,(\ref{L-3.4-pe1}), 
$f^m h\in \im D$ when $m\gg 0$. 
Then by Definition \ref{Def-MS}, we see that 
$\im D$ is indeed a Mathieu 
subspace of $k[z]$. 
\epfv

Now we can prove the main theorem of this section as follows. 
\vskip3mm
\underline{\it Proof of Theorem \ref{MainThm-1}:} 
First, by Proposition \ref{LFDer}, we only need to show that $\im D$ is a Mathieu subspace of $A$ in each of the four cases in Proposition \ref{LFDer}. Furthermore, by Lemma \ref{L-3.3} it only remains to prove case $i)$.  
So assume $D=(ax+by)\p_x+(cx+dy)\p_y$ 
for some $a, b, c, d\in k$. 

Second, by Lemma \ref{Lef-Lemma}, we may assume 
that $k$ is algebraically closed.  

Third, note that $D$ preserves the subspace 
$H\!:=kx+ky\subset A$, so its restriction 
$D|_H$ on $H$ is a linear endomorphism 
of $H$. Since $k$ is algebraically closed, 
there exists a linear automorphism $\sigma$ of 
$H$ such that the conjugation 
$\sigma (D|_H) \sigma^{-1}$ gives the Jordan 
form of $D|_H$. 
Let $\tilde \sigma$ be the unique extension 
of $\sigma$ to an automorphism of $A$. Then 
it is easy to see that 
$\tilde \sigma D \tilde \sigma^{-1}$
is also a $k$-derivation 
of $A$.

Note that 
$\im \tilde \sigma D \tilde \sigma^{-1}= \tilde\sigma(\im D)$ 
and in general Mathieu subspaces are 
preserved by $k$-algebra automorphisms. 
Therefore, we may replace $D$ by 
$\tilde \sigma D \tilde \sigma^{-1}$, 
if necessary, and assume that 
$D=a(x\p_x+y\p_y)+x\p_y$ 
(in case that the Jordan form of $D|_H$ is an $2\times 2$ 
Jordan block) or $D=ax\p_x+by\p_y$ 
(in case that the Jordan form of 
$D|_H$ is diagonal).

For the former case, by Lemma \ref{L-3.3}, $(b)$ with $m=1$, 
we see that $\im D$ is an ideal, and hence a Mathieu subspace of $A$. For the latter case, it follows from Lemma \ref{L-3.4} that $\im D$ also a Mathieu subspace of $A$. Therefore, the theorem holds.
\epfv

\renewcommand{\theequation}{\thesection.\arabic{equation}}
\renewcommand{\therema}{\thesection.\arabic{rema}}
\setcounter{equation}{0}
\setcounter{rema}{0}

\section{\bf Connection with the Two-Dimensional 
Jacobian Conjecture}\label{S4}

In the previous section we showed that 
the image of every locally finite $k$-derivation 
of $A$ is a Mathieu subspace of $A$. 
However, as we have shown in 
Example \ref{Exam2.4}, $\im D$ needs not 
to be a Mathieu subspace  
of $A$ for every $k$-derivation $D$ of $A$. This leads to the 
question of which $k$-derivations $D$ of $A$ have the property that
$\im D$ is a Mathieu subspace  of $A$. 
More precisely, we can ask

\begin{quest}\label{Quest-1}
Let $D$ be any $k$-derivation 
of $A$ such that $\dive D=0$, where for any 
$D=p\p_x+q\p_y$ $(p, q\in A)$,  
$\dive D\!:=\p_xp+\p_yq$. 
Is $\im D$ a Mathieu subspace  
of $A$?
\end{quest}

Adding one more condition, we get 

\begin{quest}\label{Quest-2}
Let $D$ be any $k$-derivation 
of $A$ such that $\dive D=0$. If $1\in \im D$, 
is $\im D$ a Mathieu subspace of $A$?
\end{quest}

Note that by Lemma \ref{OneLemma}, this question 
is equivalent to asking if $D$ is surjective 
under the further condition $1\in \im D$.

The motivation of the two questions above come from 
the following theorem.

\begin{theo}\label{MainThm2}
Question \ref{Quest-2} has an affirmative answer iff 
the two dimensional Jacobian conjecture is true.
\end{theo}

\pf $(\Rightarrow)$ Assume that Question \ref{Quest-2} has an affirmative answer. Let $F=(f, g)\in k[x, y]^2$ with 
$\det JF=1$. Consider the $k$-derivation  
$D\!:= g_y\p_x-g_x\p_y$. Then 
$\dive D=0$ and $1=\det JF=D f \in \im D$.
Since by our hypothesis $\im D$ is a Mathieu subspace  
of $A$, it follows from Lemma \ref{OneLemma} that $\im D=A$, 
i.e., $D$ is surjective. Then it follows from a 
theorem of Stein \cite{S} (see also \cite{C}) 
that $D$ is locally nilpotent. 

Since $D=\p/\p f$, $\ker D=\ker \p/\p f=k[g]$ 
by Proposition $2.2.15$ in \cite{E1}. Since $D$ has a slice $f$, 
it follows that $A=k[g][f]$, i.e., $F$ is invertible over $k$. 
So the two-dimensional Jacobian conjecture is true.

$(\Leftarrow)$  Assume that the two-dimensional Jacobian conjecture is true. Let $D=p\p_x+q\p_y$ 
$(p, q\in A)$ be a $k$-derivation of $A$ 
such that $\dive D=0$ and $1\in \im D$.

Since $\dive D=0$, we have $\p_x p=\p_y(-q)$. Then by 
Poincar\'e's lemma, there exists $g\in A$ such that 
$p=\p_y g$ and $q=-\p_x g$. 
So $D=g_y\p_x-g_x\p_y$.

Since $1\in \im D$, we get $1=D f$ for 
some $f\in A$. Let $F\!:=(f, g)\in k[x, y]^2$. 
Then we have $\det JF=D f=1$. Since by our hypothesis 
$F$ is invertible, it follows that $k[x, y]=k[f, g]$. 
Hence, we have
\begin{align*}
\im D=\im \frac{\p}{\p f}
=\frac{\p}{\p f}(k[f, g])= k[f, g]=A.
\end{align*} 
In particular, $\im D$ is a Mathieu 
subspace of $A$.
\epfv

%
%
%
%
%

\end{document}